\documentclass[10pt,a4paper]{article}
\usepackage{amsfonts}
\usepackage{amssymb}
\usepackage{mathrsfs}
\usepackage{amsthm}
\usepackage{setspace}

\usepackage{fancyhdr}
\newtheorem{defi}{Definition}[section]
\newtheorem{theo}[defi]{Theorem}
\newtheorem{prop}[defi]{Proposition}

\newtheorem{lemm}[defi]{Lemma}

\theoremstyle{definition}
\newtheorem*{rem}{Remark}

\newtheorem{exa}{Example}

\newcommand{\be}{\begin{eqnarray*}}
\newcommand{\ee}{\end{eqnarray*}}
\newcommand{\beqa}{\begin{eqnarray}}
\newcommand{\eeqa}{\end{eqnarray}}
\newcommand{\ba}{\begin{array}}
\newcommand{\ea}{\end{array}}

\begin{document}

\title{Ricci-Flat Holonomy: a Classification}
\author{Stuart Armstrong}
\date{12 September 2006}
\maketitle

\begin{abstract}
The reductive holonomy algebras for a torsion-free affine connection are analysed, with the goal of establishing which ones can correspond to a Ricci-flat connection with the same properties. Various families of holonomies are eliminated through different algebraic means, and examples are constructed (in this paper and in `Projective Geometry II: Holonomy Classification', by the same author) in the remaining cases, thus solving this problem completely, for reductive holonomy.
\end{abstract}

\tableofcontents

\section{Introduction}

Merkulov's and Schwachh\"ofer's have published the full list of reductive holonomy algebras for torsion-free, non-symmetric affine connections \cite{CIH} and \cite{CIH2}. An interesting, and hitherto unresolved question, is which of these holonomy algebras can correspond to connections which are also Ricci-flat.

In the metric case, for instance, it is well known that Levi-Civita connections with holonomies $\mathfrak{su}(p,q)$ and $\mathfrak{sp}(p,q)$ must be Ricci-flat \cite{C-Y}, whereas those with holonomy $\mathfrak{u}(p,q)$ and $\mathfrak{sp}(p,q) \oplus \mathfrak{sp}(1)$ cannot be Ricci-flat; we shall call these algebras of Ricci-type. Those with holonomy $\mathfrak{so}(n)$ may be Ricci-flat or not, and we shall attempt to extend classify all holonomy algebras into one of these three categories: Ricci-flat, Ricci-type or neither.

It is also well known that Ricci-flat symmetric spaces must have reduced holonomy \cite{Ricsy} (in the definite signature case, they must be flat, since the Ricci tensor must be a non-zero multiple of the Killing form on the Lie algebra restricted to a non-degenerate subspace, \cite{symricci2}). Thus we need only look at those irreducible holonomy algebras which are non-symmetric.

The results are summarised in tables \ref{table:17} and \ref{table:21}, the first one giving the algebras whose connections \emph{must} be Ricci-flat, and the second those who's algebras may be Ricci-flat or not.

\begin{table}[htbp]
\begin{center}
\begin{tabular}{|c|c|c||c|c|}
\hline
\hline
algebra $\mathfrak{g}$ & representation $V$ & \ \ restrictions \ \ & algebra $\mathfrak{g}$ & representation $V$ \\ 
\hline 
\hline
& & & & \\

$\mathfrak{sl}(n, \mathbb{H})$ & $\mathbb{H}^{n}$ & $ n \geq 1 $ & $\mathfrak{g}_2 (\mathbb{C})$ & $\mathbb{C}^7$ \\& & & & \\

$\mathfrak{su}(p,q)$ & $\mathbb{C}^{(p,q)}$ & $p+q \geq 3, \ (p,q)=(2,0) $ & $\mathfrak{spin}(7)$ & $\mathbb{R}^8$ \\ & & & & \\

$\mathfrak{sp}(p,q)$ & $\mathbb{H}^{(p,q)}$ & $p+q \geq 2 $ & $\mathfrak{spin}(4,3)$ & $\mathbb{R}^{(4,4)}$ \\ & & & & \\

$\mathfrak{g}_2$ & $\mathbb{R}^7$ & & $\mathfrak{spin}(7, \mathbb{C})$ & $\mathbb{C}^8$ \\ & & & & \\

$\widetilde{\mathfrak{g}}_2$ & $\mathbb{R}^{(4,3)}$ & & & \\
& & & & \\

\hline
\end{tabular}
\end{center}
\caption{Holonomy algebras that must be Ricci-flat}
\label{table:17}
\end{table}
\begin{table}[htbp]
\begin{center}
\begin{tabular}{|c|c|c|c|}
\hline
\hline
algebra $\mathfrak{g}$ & representation $V$ & \ \ restrictions \ \  \\ 
\hline 
\hline
& &  \\

$\mathfrak{so}(p,q)$ & $\mathbb{R}^{(p,q)}$ & $p+q \geq 4 $  \\ & &  \\

$\mathfrak{so}(n, \mathbb{C})$ & $\mathbb{R}^{n}$ & $ n \geq 4 $ \\ & & \\

$\mathfrak{sl}(n, \mathbb{R})$ & $\mathbb{R}^{n}$ & $n \geq 3 $  \\  & & \\

$\mathfrak{sl}(n, \mathbb{C})$ & $\mathbb{C}^{n}$ & $ n \geq 2 $ \\ & & \\

$\mathfrak{sp}(2n,\mathbb{R})$ & $\mathbb{R}^{2n}$ & $n \geq 2 $ \\ & & \\

$\mathfrak{sp}(2n,\mathbb{C})$ & $\mathbb{C}^{2n}$ & $n \geq 2 $ \\ & & \\

\hline
\end{tabular}
\end{center}
\caption{Holonomy algebras that may be Ricci-flat}
\label{table:21}
\end{table}

We shall arrive at these tables algebraically, by using various methods to exclude families of holonomy algebras. The first section will set up and define the formal curvature modules that were used by Berger and in \cite{CIH}; the Bianchi identities define the their properties and their relationship with the curvature of a torsion-free connection.

Initial results implied by the Bianchi identity then demonstrate that a Ricci-flat torsion-free connection must preserve a real volume form and, when appropriate, a complex volume form.

The next family to be dealt with were those that preserve some symplectic form; building on results from \cite{exohol} and \cite{symricci} it is not hard to show that these holonomies (apart from the maximal $\mathfrak{sp}(2n, \mathbb{R})$ and $\mathfrak{sp}(2n, \mathbb{C})$) are of Ricci-type, hence cannot correspond to Ricci-flat connections. This result is directly contained in \cite{symricci}.

To continue, we will need to build the formal curvature module of the Lie algebra $\mathfrak{g}$, the kernel of the map $V^* \otimes V^* \otimes \mathfrak{g} \to \wedge^3 V^* \otimes V$ where one anti-symmetricises over the three $V^*$. This module decomposes as $\mathfrak{g}^{(1)} \oplus H^{1,2}(\mathfrak{g})$, where $H^{1,2}(\mathfrak{g})$ is a spencer cohomology group. The module $\mathfrak{g}^{(1)}$ is the structure module for the bundle of torsion-free connections coming from the same principal $\mathfrak{g}$-bundle.

The next section deals with the various `split' algebras, those whose representation is (a sub-module of) a tensor product. We use various algebraic tricks to show that all these algebras can be dealt with in the same way, whether they be real, complex or quaternionic, symmetric, skew or hermitian. We construct a direct isomorphism between $V^* \otimes \mathfrak{g}^{(1)}$ and the Ricci curvature module. If we ignore volume forms, all split algebras, except for the minimal Segre algebras, have flat structure bundles. Consequently, they must be of Ricci-type, as must all their sub-algebras.

Then we shall look at these `minimal' Segre algebras; in other words, those whose representations are $\mathbb{R}^m \otimes_{\mathbb{R}} \mathbb{R}^2$, $\mathbb{C}^m \otimes_{\mathbb{C}} \mathbb{C}^2$ or $\mathbb{H}^m \otimes_{\mathbb{H}} \mathbb{H}^1$. Using modifications of ideas that were developed (\cite{quaternionic}) to show that a Ricci-flat $\mathbb{H}^m \otimes_{\mathbb{H}} \mathbb{H}^1$ connection must have holonomy reducing to $\mathbb{H}^m$, we similarly show that any Ricci-flat minimal Segre connection must have holonomy that reduces.

One last family remains: that of representations of $E_6$, $\mathbb{R} \oplus E_6$ and their complexifications. That second algebra also has an isomorphism between $V^* \otimes \mathfrak{g}^{(1)}$ and the Ricci tensor, making it of Ricci type. This is demonstrated by extensive algebraic manipulations. Note that although $E_6$ and $\mathbb{R} \oplus E_6$ have many formal similarities with split algebras, the methods of proof are very different.

We shall not construct explicit examples of Ricci-flat connection with the remaining holonomies. This is done mainly in papers \cite{mepro1} and \cite{mepro2}, as an adjunct to proving the existence of various projective Tractor holonomies. In fact, projective Tractor holonomies correspond to Ricci-flat Torsion-free affine holonomies on a cone one dimension higher, and \cite{mepro2} constructs these Ricci-flat cones for all higher dimensional holonomy algebras that are not ruled out by this paper.

However, there does remain the issue of some lower dimensional algebra that may correspond to Ricci-flat Torsion-free connection but not to cones. These are dealt with in the last section of this paper.

The author would like to thank Dr. Nigel Hitchin, under whose supervision and inspiration this paper was crafted. This paper appears as a section of the author's Thesis \cite{methesis}.

\section{Formal curvature modules}

\subsection{Spencer cohomology} \label{spe:co}

See \cite{CIH} for an introduction to Spencer cohomology.

Let $\mathfrak{h}$ be a Lie algebra and $V$ an $\mathfrak{h}$-module. Since $\mathfrak{h} \subset \mathfrak{gl}(V) = V \otimes V^*$, we can inductively define the modules:
\begin{eqnarray*}
\mathfrak{h}^{(-1)} &=& V \\
\mathfrak{h}^{(0)} &=& \mathfrak{h} \\
\mathfrak{h}^{(k)} &=& [ \mathfrak{h}^{(k-1)} \otimes V^* ] \cap [V \otimes \odot^{k+1} V^* ].
\end{eqnarray*}
Furthermore, if $C^{k,l}(\mathfrak{h}) = \mathfrak{h}^{(k)} \otimes \wedge^{l-1} V^*$ we may define the map
\begin{eqnarray*}
\partial : C^{k,l}(\mathfrak{h}) \to C^{k-1, l+1}(\mathfrak{h}),
\end{eqnarray*}
via anti-symmetrisation on the last $l$ indices. Since $\partial^2 = 0$, there is a complex
\begin{eqnarray*}
C^{k+1, l-1}(\mathfrak{h}) \stackrel{\partial}{\longrightarrow} C^{k,l}(\mathfrak{h}) \stackrel{\partial}{\longrightarrow} C^{k-1, l+1}(\mathfrak{h})
\end{eqnarray*}
whose cohomology at the centre term is defined to be $H^{k,l}(\mathfrak{h})$. This is called the $(k,l)$ Spencer cohomology group for $(\mathfrak{h},V)$.

\subsection{Formal curvature modules}

Given an algebra $\mathfrak{g}$ and a faithful representation $V$, there is a naturally defined operator 
\begin{eqnarray*}
\partial \left( \wedge^2 V^* \otimes \mathfrak{g} \right) \to \wedge^3 V^* \otimes V,
\end{eqnarray*}
considering $\mathfrak{g}$ as a subset of $V^* \otimes V$, $\partial$ is just antisymmetrisation over the three components. Then we define
\begin{eqnarray*}
K(\mathfrak{g}) = \mathrm{ker} \ \partial.
\end{eqnarray*}
In other words $K(\mathfrak{g})$ obeys the first Bianchi identity. The point of this construction is clear; if there is a torsion-free connection $\nu$ on a principal frame bundle $\mathcal{G}$ of the tangent bundle, then the curvature of $\nu$ is a section of
\begin{eqnarray*}
\mathcal{G} \times_G K(\mathfrak{g}).
\end{eqnarray*}
Hence we can deduce algebraic facts about the curvature of a $G$-connection from the module $K(\mathfrak{g})$.

By our results on Spencer cohomology from Section \ref{spe:co}, we know that
\begin{eqnarray*}
0 \longrightarrow \partial \left( V^* \otimes \mathfrak{g}^{(1)} \right) \longrightarrow K(\mathfrak{g}) \longrightarrow H^{1,2}(\mathfrak{g}) \longrightarrow 0.
\end{eqnarray*}
Since we will be dealing with reductive $\mathfrak{g}$'s, there actually is a splitting
\begin{eqnarray*}
K(\mathfrak{g}) = \partial \left( V^* \otimes \mathfrak{g}^{(1)} \right) \oplus H^{1,2}(\mathfrak{g}).
\end{eqnarray*}
Both of these components have a geometric interpretation; the obstruction for the given $G$-structure being flat, given that it is 1-flat -- equivalently, $M$ admitting a flat connection with principal bundle $\mathcal{G}$, given that it admits a torsion-free one -- lies in
\begin{eqnarray*}
\mathcal{G} \times_G H^{1,2}(\mathfrak{g})
\end{eqnarray*}
whereas different torsion-free connections preserving the $G$-structure differ by sections of
\begin{eqnarray*}
\mathcal{G} \times_G \mathfrak{g}^{(1)}.
\end{eqnarray*}
\begin{rem}
It is rare for an algebra to have both a $\mathfrak{g}^{(1)}$ and an $H^{1,2}(\mathfrak{g})$ component -- both an obstruction to integrability and a wide class of associated connections -- though a few do, such as the conformal $\mathbb{R}.\mathfrak{so}(p,q)$ and the almost Grassmannian $\mathbb{F}.\mathfrak{sl}(n, \mathbb{F}).\mathfrak{sl}(2, \mathbb{F})$.
\end{rem}

The full list of complex algebras with non-zero $\mathfrak{g}^{(1)}$ is as follows \cite{CIH}:
\begin{eqnarray*} \label{V:nontrivial}
\begin{array}{|c|c|c|}
\hline
\hline
\mathrm{Algebra} \ \mathfrak{g} & \mathrm{Representation} \ V & \mathfrak{g}^{(1)} \\
\hline
\hline
&& \\
\mathfrak{sl} (n, \mathbb{C}) & \mathbb{C}^n, n \geq 2 & \left( V \otimes \odot^2 V^* \right)_0 \\
& & \\
\mathfrak{gl}(n, \mathbb{C}) & \mathbb{C}^n, n \geq 1 & V \otimes \odot^2 V^* \\
& & \\
\mathfrak{gl}(n, \mathbb{C}) & \odot^2 \mathbb{C}^n, n \geq 2 & V^* \\
& & \\
\mathfrak{gl}(n, \mathbb{C}) & \wedge^2 \mathbb{C}^n, n \geq 5 & V^* \\
& & \\
\mathfrak{gl}(m, \mathbb{C}) \oplus \mathfrak{gl}(n, \mathbb{C}) & \mathbb{C}^m \otimes \mathbb{C}^n, m,n \geq 2 & V^* \\
& & \\
\mathfrak{sp}(2n, \mathbb{C}) & \mathbb{C}^{2n}, n \geq 2 & \odot^3 V^* \\
& & \\
\mathbb{C}^* \oplus \mathfrak{sp}(2n, \mathbb{C}) & \mathbb{C}^{2n}, n \geq 2 & \odot^3 V^* \\
& & \\
\mathfrak{co}(n, \mathbb{C}) & \mathbb{C}^n, n \geq 5 & V^* \\
& & \\
\mathbb{C}^* \oplus \mathfrak{spin}(10, \mathbb{C}) & \mathbb{C}^{16} & V^* \\
& & \\
\mathbb{C}^* \oplus \mathfrak{e}_6^{\mathbb{C}} & \mathbb{C}^{27} & V^* \\
\hline
\end{array}
\end{eqnarray*}

Since $K(\mathfrak{g})$ is a formal curvature module, we may define the formal Ricci-curvature module $R(\mathfrak{g})$ by taking the trace of $K(\mathfrak{g})$. Then possible Ricci-flat curvatures will lie inside the kernel of
\begin{eqnarray*}
K(\mathfrak{g}) \to R(\mathfrak{g}).
\end{eqnarray*}

If, on the other hand, this map has no kernel; in other words, if
\begin{eqnarray*}
K(\mathfrak{g}) \cong R(\mathfrak{g}),
\end{eqnarray*}
then we say that $\mathfrak{g}$ has curvature of Ricci-type. Obviously a connection whose holonomy algebra is of Ricci-type cannot be Ricci-flat without being flat.

\subsection{Complex modules} \label{complex:algebras}

Let $(V, J)$ be a vector space with complex structure $J$. By $\overline{V}$ we mean $(V, -J)$. Let $\alpha : \otimes^j V \to \wedge^j V$ be the natural antisymmetrisation map. Let $\otimes^{(n,m)} V = \otimes_{\mathbb{C}}^n V \otimes_{\mathbb{C}} \otimes_{\mathbb{C}}^m \overline{V}$. Then we shall define the space $\wedge^{(n,m)}$ as:
\begin{defi}
\begin{eqnarray*}
\wedge^{(n,m)} V = \alpha \left( \otimes^{(n,m)} V \right).
\end{eqnarray*}
\end{defi}
Obviously, $\wedge^{(n,m)} V = \overline{\wedge^{(m,n)} V}$, implying that $\wedge^{(n,m)} V$ and $\wedge^{(m,n)} V$ are the same spaces. And, of course, $\wedge^{(n,0)} V = \wedge^n_{\mathbb{C}} V$. It can easily be seen that
\begin{lemm} \label{little:lemma}
$\partial (\wedge^{(n,m)} V \otimes V) \subset \wedge^{(n+1,m)} V \oplus \wedge^{(n,m+1)} V$ if $n \neq m$, and $\partial (\wedge^{(n,n)} V \otimes V) \subset \wedge^{(n+1,1)} V$.
\end{lemm}

For the rest of this section, any tensor product is complex unless stated otherwise. Let $\mathfrak{g}_{\mathbb{R}}$ be a real Lie algebra, with a corresponding complex Lie algebra $\mathfrak{g}$. Let $V_{\mathbb{R}}$ be a real representation of $\mathfrak{g}_{\mathbb{R}}$, and $V = V_{\mathbb{R}} \otimes_{\mathbb{R}} \mathbb{C}$ the corresponding representation of $\mathfrak{g}$. For any two complex spaces $W$ and $U$,
\begin{eqnarray*}
W \otimes_{\mathbb{R}} U = (W \otimes U) \ \oplus \ (\overline{W} \otimes U),
\end{eqnarray*}
the $+1$ and $-1$ eigen-spaces of the operator $J \otimes J$. Similarly, the module $\wedge^2_{\mathbb{R}} V^* \otimes_{\mathbb{R}} \mathfrak{g}$ splits into three sub-modules:
\begin{eqnarray*}
\wedge^2_{\mathbb{R}} V^* \otimes_{\mathbb{R}} \mathfrak{g} &=&
\left( \wedge^{(2,0)} V^* \otimes \mathfrak{g} \right) \oplus
\left( \wedge^{(1,1)} V^* \otimes_{\mathbb{R}} \mathfrak{g} \right) \oplus
\left( \wedge^{(0,2)} V^* \otimes \mathfrak{g} \right)
\\ &=& P_1 \oplus P_2 \oplus P_3,
\end{eqnarray*}
where $\wedge^{(2,0)} V^* = \wedge^2_{\mathbb{C}} V^*$ and $\wedge^{(0,2)} V^*$ is the same space with opposite complex structure. The space $\wedge^{(1,1)} V^*$ is just the space of skew-hermitian forms; this space does \emph{not} have a complex structure itself, hence the real tensor product in the central term. Denote by $p_1, p_2, p_3$ the projections onto these sub-modules. These modules are disjoint from the point of view of the $\partial$ map:
\begin{lemm} \label{split:delta}
If $\partial (a) = 0$, then $\partial p_j(a) = 0$ for all $j$.
\end{lemm}
\begin{proof}
Assume $\partial (a) = 0$. The module $P_1$ is contained in the module $\otimes^3 V^* \otimes V$, so $\partial (P_1) \subset \wedge^{(3,0)} V^* \otimes V$. By Lemma \ref{little:lemma}, $\partial(P_2)$ and $\partial (P_3)$ are both contained in $\wedge^{2,1} V^* \otimes_{\mathbb{R}} V$. Consequently, $\partial p_1(a)$ must be zero. From now on, by replacing $a$ with $a - p_1(a)$, we may assume that $p_1(a) = 0$.

The operator $\Theta = J \otimes J \otimes J \otimes J$ operates naturally on $\wedge^2_{\mathbb{R}} V^* \otimes_{\mathbb{R}} V^* \otimes_{\mathbb{R}} V$, and, since $\partial$ is an antisymmetrisation of this space and $\Theta$ is entirely symmetric,
\begin{eqnarray*}
\Theta \circ \partial = \partial \circ \Theta.
\end{eqnarray*}
However, $\Theta(p_2(a)) = -p_2(a)$ and $\Theta(p_3(a)) = p_3(a)$, so
\begin{eqnarray*}
&\partial p_2(a) = \frac{1}{2} \left( \partial(a) - \Theta \partial(a) \right) = 0& \\ &\mathrm{and}& \\
&\partial p_3(a) = \frac{1}{2} \left( \partial(a) + \Theta \partial(a) \right) = 0.&
\end{eqnarray*}
\end{proof}
On the other hand, $\wedge^{(2,0)} V^* \otimes \mathfrak{g}$ is just the complexification of the real module $\wedge^2_{\mathbb{R}} V^*_{\mathbb{R}} \otimes_{\mathbb{R}} \mathfrak{g}_{\mathbb{R}}$. So we can directly classify this piece of the complex module in terms of the real one:
\begin{prop} $ p_1 \big( K(\mathfrak{g})  \big) = K(\mathfrak{g}_{\mathbb{R}}) \otimes_{\mathbb{R}} \mathbb{C}$. \end{prop}
The next lemma deals with the $P_3$ component:
\begin{lemm} $\partial$ is injective on $P_3$. \end{lemm}
\begin{proof}
Let $b_1$ be an element of $P_3$. Then $\partial (b_1)$ equals $\frac{1}{3} (b_1 + b_2 + b_3)$ where $b_2$ and $b_3$ are the two cyclic permutations of $b_1$. However, if we apply $\theta = J \otimes J$ to the first two components of these elements, we see that:
\begin{eqnarray*}
\theta b_1 &=& -b_1 \\
\theta b_2 &=& b_2 \\
\theta b_3 &=& b_3.
\end{eqnarray*}
Accordingly, $b_1 = \frac{3}{2} \big( \partial(b_1) - \theta \partial(b_1) \big)$, directly displaying the injectivity of $\partial$ on $P_3$.
\end{proof}
Putting this together with Lemma $\ref{split:delta}$ implies that $p_3 (a)$ must be zero if $\partial(a) = 0$. In other words,
\begin{eqnarray*}
p_3 \big( K( \mathfrak{g}) \big) = 0.
\end{eqnarray*}
Thus:
\begin{theo} \label{split:theo}
The formal curvature module $K(\mathfrak{g})$ splits as
\begin{eqnarray*}
K(\mathfrak{g}) &=& K_1(\mathfrak{g}) \oplus K_2(\mathfrak{g}),
\end{eqnarray*}
where $K_1(\mathfrak{g})$ is the complexification of $K(\mathfrak{g}_{\mathbb{R}})$ and $K_2(\mathfrak{g}) \subset \wedge^{(1,1)} V^* \otimes_{\mathbb{R}} \mathfrak{g}$.
\end{theo}
Furthermore, the formal Ricci module splits into the sum of the traces of these two modules:
\begin{eqnarray*}
R(\mathfrak{g}) &=& R_1(\mathfrak{g}) \oplus R_2(\mathfrak{g}),
\end{eqnarray*}
with $R_1(\mathfrak{g})$ a $J$-symmetric space, and $R_2(\mathfrak{g})$ a $J$-hermitian space.

Note that since this splitting result is true for $\mathfrak{gl}(n, \mathbb{C})$, it is also true for any $\mathfrak{g} \subset \mathfrak{gl}(n, \mathbb{C})$, even if $\mathfrak{g}$ is not itself a complex algebra (such as $\mathfrak{u}(n)$).

\begin{exa} \label{so:nc}
To illustrate these two bundles, we can use two metric examples; first of all, let $\mathfrak{g} = \mathfrak{so}(n,\mathbb{C})$. The complex metric gives an isomorphism $\mathfrak{g} \cong \wedge^{(2,0)} V^*$, and the extra metric condition that $R_{hjkl} = R_{klhj}$ gives us
\begin{eqnarray*}
K_1(\mathfrak{g}) &\subset & \mathfrak{g} \otimes \mathfrak{g}, \\
K_2(\mathfrak{g}) &=& 0.
\end{eqnarray*}
And, of course, the Ricci tensor of such a connection must be $J$-symmetric.
\end{exa}

\begin{exa} \label{u:n}
Conversely, let $\mathfrak{g} = \mathfrak{u}(n)$. The hermitian metric gives an isomorphism $\mathfrak{g} = \wedge^{(1,1)} V^*$, and with the condition $R_{hjkl} = R_{klhj}$ as before,
\begin{eqnarray*}
K_1(\mathfrak{g}) &=& 0, \\
K_2(\mathfrak{g}) &\subset & \mathfrak{g} \otimes_{\mathbb{R}} \mathfrak{g}.
\end{eqnarray*}
And, of course, these K\"ahler manifolds must have $J$-hermitian Ricci tensor.
\end{exa}

\section{Volume forms and the Ricci tensor}

Let $\mathcal{E}_n = \wedge^n T^*$ be the volume bundle on a manifold $M^n$, and $\nabla$ a torsion free-connection on $M$. Then the curvature $R_{hj \phantom{k} l}^{\phantom{hj}k}$ of $\nabla$ acts on $\mathcal{E}_n$ via its trace $R_{hj \phantom{k} k}^{\phantom{hj}k}$. However, since $\nabla$ is torsion-free, the first Bianchi identity gives
\begin{eqnarray*}
R_{hj \phantom{k} k}^{\phantom{hj}k} &=& R_{kj \phantom{k} h}^{\phantom{hj}k} + R_{hk \phantom{k} j}^{\phantom{hj}k} \\ &=& \mathsf{Ric}_{hj} - \mathsf{Ric}_{jh}.
\end{eqnarray*}
This demonstrates the next lemma:
\begin{lemm}
A torsion-free connection $\nabla$ preserves a volume form if and only if its Ricci tensor is symmetric.
\end{lemm}

Similarly, if $n = 2m$ and $\nabla$ preserves a complex structure, let $\mathcal{E}^{\mathbb{C}}_{m} = \wedge^{m,0} T^*_{\mathbb{C}}$ be the complex volume bundle. Then the curvature of $\nabla$ acts on $\mathcal{E}^{\mathbb{C}}_{m}$ via the complex trace 
\begin{eqnarray*}
\mathrm{trace}_{\mathbb{C}} \ R &=& \frac{1}{2} \left( \mathrm{trace}_{\mathbb{R}} R + i \ \mathrm{trace}_{\mathbb{R}} JR \right).
\end{eqnarray*}
The first term is just the skew-symmetric part of the Ricci tensor, as before. The second term is given by
\begin{eqnarray*}
R_{hj \phantom{k} l}^{\phantom{hj} k} J^l_k &=& \left( R_{lj \phantom{k} h}^{\phantom{hj} k} + R_{hl \phantom{k} j}^{\phantom{hj} k} \right) J^l_k.
\end{eqnarray*}
Since $\nabla$ preserves the complex structure, $R_{hj \phantom{k} l}^{\phantom{hj}k} J^n_k = R_{hj \phantom{k} k}^{\phantom{hj}n} J^k_l$, implying that the previous formula becomes:
\begin{eqnarray*}
R_{hj \phantom{k} l}^{\phantom{hj} k} J^l_k &=& R_{lj \phantom{k} k}^{\phantom{hj} l} J^k_h + R_{hl \phantom{k} k}^{\phantom{hj} l} J^k_j \\ &=& - \mathsf{Ric}_{jk}J^k_h + \mathsf{Ric}_{hk} J^k_j,
\end{eqnarray*}
the skew-symmetric part of $\mathsf{Ric}J$. This gives us the result:
\begin{lemm} \label{complex:ricci}
A torsion-free connection $\nabla$ preserves a complex volume form if and only if the tensors $\mathsf{Ric}$ and $\mathsf{Ric}J$ are both symmetric.
\end{lemm}

And this gives us our first tool for classifying Ricci-Flat spaces, notably that
\begin{prop}
A Ricci-flat space $(M, \nabla)$ with $\nabla$ torsion-free, has a preserved real volume form, and, if $\nabla$ preserves a complex structure, it also has a preserved complex volume form.
\end{prop}
\begin{exa}
Looking back at Example \ref{so:nc}, $\mathfrak{g} = \mathfrak{so}(n,\mathbb{C})$, we see that its Ricci tensor is $J$-symmetric. Being a metric connection, its Ricci tensor must also be symmetric, so we come to the unsurprising conclusion that a connection with holonomy $\mathfrak{so}(n,\mathbb{C})$ must preserve a complex volume form.
\end{exa}
\begin{exa}
On the other hand, Example \ref{u:n} shows that $\mathfrak{g} = \mathfrak{u}(n)$ has a Ricci tensor that is $J$-hermitian, in other words $J$-skew. This gives us the slightly more interesting conclusion that a K\"ahler manifold has a preserved complex volume form (i.e. has $\mathfrak{su}(n)$ holonomy) if and only if it is Ricci-flat.
\end{exa}
\section{Symplectic sub-algebras}
These are the various sub-algebras of the symplectic and complex symplectic algebras, $\mathfrak{sl}(2n, \mathbb{R})$ and $\mathfrak{sl}(2n, \mathbb{C})$. The list of such algebras that can appear as irreducible holonomy algebras is as follows \cite{CIH}:
\begin{eqnarray*}
\begin{array}{|c|c||c|c|}
\hline
\hline
\mathrm{Algebra} \ \mathfrak{g} & \mathrm{Representation} \ V &\mathrm{Algebra} \ \mathfrak{g} &\mathrm{Representation} \ V \\
\hline
\hline
\mathfrak{sp}(2n,\mathbb{R}) & \mathbb{R}^{2n} & \mathfrak{e}_7^5 &\mathbb{R}^{56} \\
& & & \\
\mathfrak{sp}(2n,\mathbb{C}) & \mathbb{C}^{2n} & \mathfrak{e}_7^7 &\mathbb{R}^{56} \\
& & & \\
\mathfrak{sl}(2,\mathbb{R}) & \mathbb{R}^{4}=\odot^3 \mathbb{R}^2 & \mathfrak{e}_7^{\mathbb{C}} &\mathbb{C}^{56} \\
& & & \\
\mathfrak{sl}(2,\mathbb{C}) & \mathbb{C}^{4}=\odot^3 \mathbb{C}^2 & \mathfrak{spin}(2,10) &\mathbb{R}^{32} \\ 
& & & \\
\mathfrak{sl}(2,\mathbb{R}) \oplus \mathfrak{so}(p,q) & \mathbb{R}^{2(p+q)}, p+q \geq 3& \mathfrak{spin}(6,6) &\mathbb{R}^{32} \\ 
& & & \\
\mathfrak{sl}(2,\mathbb{C}) \oplus \mathfrak{so}(n, \mathbb{C}) & \mathbb{C}^{2n}, n \geq 3 & \mathfrak{spin}(12,\mathbb{C}) &\mathbb{C}^{32} \\ 
& & & \\
\mathfrak{sp}(1) \oplus \mathfrak{so}(n, \mathbb{H}) & \mathbb{H}^{n}, n \geq 2 & \mathfrak{sp}(6,\mathbb{R}) &\mathbb{R}^{14} \subset \wedge^3 \mathbb{R}^6 \\
& & & \\
\mathfrak{sl}(6, \mathbb{R}) & \mathbb{R}^{20} \cong \wedge^3 \mathbb{R}^6 & \mathfrak{sp}(6,\mathbb{C}) &\mathbb{C}^{14} \subset \wedge^3 \mathbb{C}^6 \\ 
& & & \\
\mathfrak{sl}(6,\mathbb{C}) & \mathbb{C}^{20} \cong \wedge^3 \mathbb{C}^6 & & \\ 
& & & \\
\mathfrak{su}(1,5) & \mathbb{R}^{20} & & \\ 
& & & \\
\mathfrak{su}(3,3) & \mathbb{R}^{20} & & \\ 
\hline
\end{array}
\end{eqnarray*}

This section aims to prove the following theorem:
\begin{theo}
All the algebras in that list, apart from $\mathfrak{sp}(2n, \mathbb{R})$ and $\mathfrak{sp}(2n, \mathbb{C})$ themselves, have curvature of Ricci-type:
\begin{eqnarray*}
K(\mathfrak{g}) \cong R(\mathfrak{g}).
\end{eqnarray*}
In other words, connections with these holonomies cannot be Ricci-flat without being flat.
\end{theo}
Fix a given algebra $\mathfrak{g}$, a proper subset of $\mathfrak{sp}(V, \mathbb{F})$, $V \cong \mathbb{R}^{2n}$. There are canonical manifolds with full holonomy $\mathfrak{g}$; they are constructed in \cite{exohol} using perturbed Poisson structures, and locally any manifold with $\mathfrak{g}$-holonomy is constructed in this way. However, we shall not need this explicit construction, as we shall demonstrate this theorem algebraically.

Fix a given symplectic form $\eta \in \wedge^2 V^*$. Given $\eta$, and since $\mathfrak{g}$ is semi-simple, we have a $\mathfrak{g}$-invariant projection
\begin{eqnarray*}
\odot^2 V \to \mathfrak{g}.
\end{eqnarray*}
Call $u \circ v$ the projection of $u \odot v$. Then, by \cite{exohol}, \cite{symricci} and \cite{E7}, the following equalities hold for all $\mathfrak{g}$ in the list:
\begin{eqnarray}
\label{sympl:one} \eta(Au,v) &=& (A, u \circ v) \\
\label{sympl:two} (u \circ v , s \circ t) - (u \circ t, s \circ v) &=& \big( 2 \eta(u,s) \eta(v,t) \\
&& + \nonumber \eta(u,t) \eta(v,s) + \eta(u,v)\eta(s,t) \big),
\end{eqnarray}
for all $A \in \mathfrak{g}$ and all $u,v,s,t \in V$, with
\begin{eqnarray*}
(-,-) = -\frac{1}{4n+4} B
\end{eqnarray*}
where $B$ is the Killing form on $\mathfrak{g}$ (which is the restriction of the Killing form on $\mathfrak{sp}(V, \mathbb{F})$). There is an injection of $Ad(\mathfrak{g})$ into $K(\mathfrak{g})$ given by $A \to \rho_A$,
\begin{eqnarray*}
\begin{array}{rccl}
\rho_A: & \wedge^2 V &\longrightarrow & \mathfrak{g} \\
& u \wedge v & \longrightarrow & 2 \eta(u,v)A - u \circ (Av) + v \circ (Au).
\end{array}
\end{eqnarray*}
The fact that $\rho_A \in K(\mathfrak{g})$ is guaranteed by Equations \ref{sympl:one} and \ref{sympl:two}. Paper \cite{CIH} demonstrates that the whole of $K(\mathfrak{g})$ in constructed in this manner. Then we have the following Proposition, coming from \cite{symricci}:
\begin{prop}
$\mathsf{Ric}(\rho_A) = 0$ iff $A= 0$.
\end{prop}
\begin{proof}
We shall use the following lemma:
\begin{lemm}
For any element $k \in K(\mathfrak{sp}(V, \mathbb{F}))$,
\begin{eqnarray*}
\mathsf{Ric}(k)(x,y) = \eta (k(\eta^{-1})x,y).
\end{eqnarray*}
\end{lemm}
\begin{proof}
Let $(e_j, f_j)$ be a basis for $V$ such that, when using the summation convention, $\eta^{-1} = e_j \wedge f_j$. Thus, continuing with the summation convention,
\begin{eqnarray*}
\mathsf{Ric}(k)(x,y) &=& tr (k(x,-)y) = \eta(R_{(e_j ,x)}y, f_j) - \eta(R_{(f_j,x)}y, e_j) \\
&=& - \eta(R_{(x, e_j )}f_j, y) - \eta( R_{(f_j ,x)} e_j,y) \\
&=& \eta(R_{(e_j, f_j)} x, y),
\end{eqnarray*}
as $\eta$ maps $k$ to an element of $\wedge^2 V^* \otimes \odot^2 V^*$.
\end{proof}
Now suppose $\mathsf{Ric}(\rho_A) = 0$. This is the case iff $\rho_A (\eta^{-1}) = 0$. But then \cite{symricci} demonstrates $\rho_A (\eta^{-1}) = 0$ only when $A = 0$.
\end{proof}
We have consequently shown that $K(\mathfrak{g}) \cong R(\mathfrak{g})$, or in other words that $\mathfrak{g}$ is of Ricci-type.

\section{Split spaces: General case}
Let $(M, \nabla)$, be a manifold with affine connection, whose holonomy algebra bundle acts irreducibly on $T$. Let $\mathfrak{g}$ be the fiber of the holonomy algebra at a point, and $V$ the fiber of $T$ at the same point. By our assumptions, $V$ is an irreducible representation of $\mathfrak{g}$.

Then we call $M$ a split space if $V$ is in some way the tensor product of smaller representations of $\mathfrak{g}$. In details, we say that $\mathfrak{g}$ is a \emph{maximal} algebra if there does not exist a non-symmetric holonomy algebra $\mathfrak{h}$ such that $\mathfrak{g}$ is a strict subalgebra of $\mathfrak{h}$ and
\begin{eqnarray*}
[\mathfrak{h}, \mathfrak{h}] = [\mathfrak{g}, \mathfrak{g}].
\end{eqnarray*}
More intuitively, $\mathfrak{g}$ is maximal if it has the maximal allowed reductive piece. For instance, $\mathfrak{gl}(n)$, $\mathfrak{co}(n)$ and $\mathfrak{u}(n)$ are maximal, whereas $\mathfrak{sl}(n)$, $\mathfrak{so}(n)$ and $\mathfrak{su}(n)$ are not. The algebra $\mathfrak{spin}(7)$ is also maximal, since $\mathbb{R} \oplus \mathfrak{spin}(7)$ is not a possible holonomy algebra.

Then the following table gives the maximal split algebras:

\begin{eqnarray*}
\begin{array}{|c|c|c|}
\hline
\hline
\mathrm{Algebra} \ \mathfrak{g} & \mathrm{Representation} \ V & \mathrm{Restrictions} \\
\hline
\hline
& & \\
\mathfrak{gl}(n, \mathbb{R}) & \odot^2 \mathbb{R}^n & n \geq 3  \\
& & \\
\mathfrak{gl}(n, \mathbb{R}) & \wedge^2 \mathbb{R}^n & n \geq 5  \\
& & \\
\mathfrak{gl}(n, \mathbb{C}) & \odot^2 \mathbb{C}^n & n \geq 3  \\
& & \\
\mathfrak{gl}(n, \mathbb{C}) & \wedge^2 \mathbb{C}^n & n \geq 5  \\
& & \\
\mathfrak{gl}(n, \mathbb{C}) & \ \ H^+_n(\mathbb{C}) \cong \wedge^{1,1} \mathbb{C}^n \ \ & n \geq 3  \\
& & \\
\mathfrak{gl}(n, \mathbb{H}) & H^+_n(\mathbb{H}) & n \geq 3  \\
& & \\
\mathfrak{gl}(n, \mathbb{H}) & H^-_n(\mathbb{H}) & n \geq 2  \\
& & \\
\mathbb{C} \oplus \mathfrak{sl}(m, \mathbb{C}) \oplus \mathfrak{sl}(r, \mathbb{C}) & \mathbb{C}^m \otimes \mathbb{C}^r & m > r \geq 2 \ \mathrm{or} \ m \geq r > 2 \\
& & \\
\mathbb{R} \oplus \mathfrak{sl}(m, \mathbb{R}) \oplus \mathfrak{sl}(r, \mathbb{R}) & \mathbb{R}^m \otimes \mathbb{R}^r & m > r \geq 2 \ \mathrm{or} \ m \geq r > 2 \\
& & \\
\mathbb{R} \oplus \mathfrak{sl}(m, \mathbb{H}) \oplus \mathfrak{sl}(r, \mathbb{H}) & \mathbb{H}^m \otimes \mathbb{H}^r \cong \mathbb{R}^{4mr} & m > r \geq 1 \ \mathrm{or} \ m \geq r > 1 \\
\hline
\end{array}
\end{eqnarray*}
Here $H_n^+(\mathbb{F})$ is the space of self-adjoint $n$ by $n$ matrices with entries in $\mathbb{F}$, whereas $H_n^-(\mathbb{F})$ is the complementary space of skew adjoint ones. Notice that under multiplication by $i$,
\begin{eqnarray*}
H_n^+(\mathbb{F}) \cong H_n^-(\mathbb{C}) = \wedge^{(1,1)} \mathbb{C}^n
\end{eqnarray*}
where $\wedge^{(1,1)} \mathbb{C}^n$ is defined as in section \ref{complex:algebras}.

All the algebras on this table share the property that
\begin{eqnarray*}
\mathfrak{g}^{(1)} &=& V^*,
\end{eqnarray*}
see Table \ref{V:nontrivial} and \cite{CIH}. Then we aim to prove the following theorem:
\begin{theo}
\begin{eqnarray*}
\partial \left( V^* \otimes_{\mathbb{R}} \mathfrak{g}^{(1)} \right) \cong R(\mathfrak{g}).
\end{eqnarray*}
\end{theo}
This is enough to specify all of the algebras on this table except for the minimal Segre ones:
\begin{eqnarray*}
\mathbb{C} \oplus \mathfrak{sl}(m, \mathbb{C}) \oplus \mathfrak{sl}(2,\mathbb{C}), \\
\mathbb{R} \oplus \mathfrak{sl}(m, \mathbb{R}) \oplus \mathfrak{sl}(2,\mathbb{R}), \\
\mathbb{R} \oplus \mathfrak{sl}(m, \mathbb{H}) \oplus \mathfrak{sl}(1,\mathbb{H}).
\end{eqnarray*}
For in all other cases the obstruction tensor $H^{1,2}(\mathfrak{g}) = 0$, \cite{CEEH}, so
\begin{theo}
All algebras on the table except for the minimal Segre are of Ricci-type. Consequently, neither they nor any of their subalgebras may be holonomy algebras of Ricci-flat connections.
\end{theo}
That theorem will remove most of what's left of possible cone holonomies.

For the rest of this section, any unspecified tensor product $\otimes$ is taken to be a real tensor product. Let $W \cong \mathbb{R}^{m}$ and $U \cong \mathbb{R}^{r}$, and let $E = \mathbb{R}^{m} \otimes \mathbb{R}^r$. Choose $(X_k)$ and $(Y_j)$, basis of $\mathbb{R}^{m}$ and $\mathbb{R}^{r}$, with dual basis $(x^k)$ and $(y^j)$. Then define $\mu: E^* \to E^* \odot E^* \otimes E$,
\begin{eqnarray*}
\mu (ab) &=& ay^j \otimes x^kb \otimes X_k Y_j + x^kb \otimes ay^j \otimes X_k Y_j,
\end{eqnarray*}
summing over repeated indexes.
\begin{lemm}
The function $\mu$ is independent of the choice of basis $(X_k)$ and $(Y_j)$, and is injective.
\end{lemm}
\begin{proof}
$\mu$ is the sum of two elements, each a reordering of the tensor product
\begin{eqnarray*}
ab \otimes x^j X_j \otimes y^k Y_k = ab \otimes Id_{\mathbb{R}^m} \otimes Id_{\mathbb{R}^r}
\end{eqnarray*}
and that element is obviously independent of the basis. For injectivity, note that the trace of $\mu(ab)$ over the last two elements is
\begin{eqnarray*}
\mathrm{trace} \ \mu (ab) &=& (m+r) ab.
\end{eqnarray*}
\end{proof}

Given complex structures $J^U$ and $J^W$ on $U$ and $W$, we can define the inclusion of the complex tensor product into the real one, $U \otimes_J W \subset E$, with $J= (J^U, J^W)$. This is the subbundle spanned by elements of the form
\begin{eqnarray*}
a \otimes b - J^Wa \otimes J^Ub.
\end{eqnarray*}
Similarly, if $U$ is a right quaternionic structure, $J_1^U J_2^U = -J_3^U$, and $W$ a left quaternionic structure, $J_1^W J_2^W = J_3^W$, we may define the quaternionic tensor product bundle $U \otimes_{\mathbb{H}} W$ inside $E$ as the intersection
\begin{eqnarray*}
\big( U \otimes_{J_1} W \big) \cap \big( U \otimes_{J_3} W \big) \cap \big( U \otimes_{J_3} W \big),
\end{eqnarray*}
$J_k = (J^U_k, J^W_k)$ as before. In fact, we need only take the intersections of the first two bundles.

Similarly, in the case when $m=r$, $W = U$, we may define the alternating $W \wedge W$ and symmetric spaces $W \odot W$ in the usual way. Then all of our representation spaces $V$ are intersections of these various bundles; for instance in the real Segre case
\begin{eqnarray*}
V = E,
\end{eqnarray*}
whereas in the complex symmetric case,
\begin{eqnarray*}
V = \big( W \otimes_{(+J,+J)} W \big) \cap \big( W \odot W \big)
\end{eqnarray*}
while the complex self-adjoint bundle is given by
\begin{eqnarray*}
V = \big( W \otimes_{(+J,-J)} W \big) \cap \big( W \odot W \big),
\end{eqnarray*}
and so on. As any complex structure $J$ has a dual action on the dual bundle, and the transpose operation applies naturally to a space and its dual, for any of our spaces $V \subset E$, we have a well defined $V^* \subset E^*$. Consequently, we have a well defined projection $p: E^* \to V^*$.

\begin{prop}
For the representation space $V^* \subset E^*$ of a split algebra $\mathfrak{g}$,
\begin{eqnarray*}
p \circ \mu (V^*) = \mathfrak{g}^{(1)},
\end{eqnarray*}
where $p$ is operating on the first element of $\mu(V^*)$.
\end{prop}
\begin{proof}
Since we know that $\mathfrak{g}^{(1)} \cong V^*$ and that $\mu$ is injective, it suffices to show that $p \circ \mu (V^*) \subset \mathfrak{g}^{(1)}$. First we shall use the lemma:
\begin{lemm}
$\mu(V^*)\llcorner (V)  \subset \mathfrak{g}$.
\end{lemm}
\begin{proof}
In the Segre case,
\begin{eqnarray*}
\mu(ab) \llcorner CD = a(C) x^kb \otimes X_k D + b(D) ay^j \otimes C Y^j.
\end{eqnarray*}
This corresponds to an element of $\mathfrak{gl}(m, \mathbb{R}) \oplus \mathfrak{gl}(r, \mathbb{R})$.

In the skew case,
\begin{eqnarray*}
\mu(ab - ba) \llcorner (CD - DC) &=& +a(C) x^kb \otimes X_k D + b(D) ay^j \otimes C Y^j
\\ && - b(C) x^ka \otimes X_k D - a(D) by^j \otimes C Y^j
\\ && - a(D) x^kb \otimes X_k C - b(C) ay^j \otimes D Y^j
\\ && + b(D) x^ka \otimes X_k C + a(C) by^j \otimes D Y^j,
\end{eqnarray*}
which corresponds to
\begin{eqnarray*}
a(C) b \otimes D - b(C) a \otimes D - a(D) b \otimes C + b(D) a \otimes C \in \mathfrak{g},
\end{eqnarray*}
acting diagonally inside $\mathfrak{gl}(m, \mathbb{R}) \oplus \mathfrak{gl}(m, \mathbb{R})$. The proof in the symmetric case is the same, modulo a few sign differences.

In the complex case:
\begin{eqnarray*}
& A = \mu(ab - JaJb) \llcorner (CD - JC JD)& \\
&=& \\
& +a(C) x^kb \otimes X_k D + b(D) ay^j \otimes C Y^j& \\
& - Ja(C) x^kJb \otimes X_k D - Jb(D) Jay^j \otimes C Y^j& \\
& - a(JC) x^k b \otimes X_k JD - b(JD) ay^j \otimes JC Y^j& \\
& + Ja(JC) x^kJb \otimes X_k JD + Jb(JD) aby^j \otimes JC Y^j&,
\end{eqnarray*}
then, using the fact that $Ja(C) = a(JC)$ and $Id = X^k \otimes x_k = - JX^k \otimes Jx_k$, one has
\begin{eqnarray*}
A (+ef) = A (-JeJf),
\end{eqnarray*}
implying that $A$ is contained in
\begin{eqnarray*}
\mathfrak{gl}(m,J) \oplus \mathfrak{gl}(r,J).
\end{eqnarray*}
The lemma in the general case then follows from intersections of these various constructions.
\end{proof}
This lemma establishes that $p \circ \mu(V^*) \subset V^* \otimes \mathfrak{g}$. Moreover, if $v \in V$, $w \in V^*$,
\begin{eqnarray*}
p \circ \mu(w)(v) = \mu(w)(v)
\end{eqnarray*}
by definition of what $p$ is. Consequently $p \circ \mu(w)$ remains symmetric in the first two elements; consequently
\begin{eqnarray*}
p \circ \mu(V^*) \subset \big(V^* \odot V^* \otimes V \big) \cap \big( V^* \otimes \mathfrak{g} \big) = \mathfrak{g}^{(1)}.
\end{eqnarray*}
\end{proof}
Before continuing, we shall see what properties $\mu(V^*)$ and $p \circ \mu(V^*)$ share; for the first is easier to work with. First of all, we know that for $v \in V$, $w \in V^*$,
\begin{eqnarray*}
p \circ \mu(w)(v) = \mu(w)(v).
\end{eqnarray*}
However, we shall also need:
\begin{lemm}
Both $\mu(V^*)$ and $p \circ \mu(V^*)$ have the same trace over the last two elements -- equivalently over the first and last element.
\end{lemm}
\begin{proof}
The projection $p$ commutes with the operation of taking traces. However, the trace formula is
\begin{eqnarray*}
\mathrm{trace} \ \mu (ab) &=& (m+r) ab.
\end{eqnarray*}
So $\mathrm{trace} \ \mu (V^*) \subset V^*$. Therefore, as $p$ is the identity on $V^*$,
\begin{eqnarray*}
\mathrm{trace} \ \mu (V^*) = p \circ \mathrm{trace} \ \mu (V^*) = \mathrm{trace} \ p \circ \mu (V^*).
\end{eqnarray*}
\end{proof}

We are now in a position to prove the main theorem. Let $\mathbf{R}$ be the operator taking the Ricci-trace. Recall that:
\begin{eqnarray*}
\partial (cd \otimes \mu(ab)) &=& (cd \otimes ay^j) \otimes x^kb \otimes X_k Y_j + (cd \otimes x^kb) \otimes ay^j \otimes X_k Y_j \\
&& - (ay^j \otimes cd) \otimes x^kb \otimes X_k Y_j - (x^kb \otimes cd) \otimes ay^j \otimes X_k Y_j
\end{eqnarray*}
\begin{lemm}
The linear maps $\mathbf{R} \ \partial (Id_{E^*} \otimes p \circ \mu )$ and $\mathbf{R} \ \partial (Id_{E^*} \otimes \mu)$, both mapping $E^* \otimes E^*$ to itself, are equal on $V^* \otimes V^* \subset E^* \otimes E^*$.
\end{lemm}
\begin{proof}
This is a direct consequence of the two identities for $p \circ \mu$ and $\mu$ that we have just seen. The first two terms of $\mathbf{R} \ \partial (V^* \otimes \mu(V^*))$ involve evaluating an element of $\mu(V^*)$ on an element of $V^*$; the second two terms involve the tensor product of an element of $V^*$ with the trace of an element of $\mu(V^*)$. And one can replace $\mu$ with $p \circ \mu$ in all these cases.
\end{proof}
Then the final statement is a consequence of:
\begin{prop}
The linear map $\mathbf{P} = \mathbf{R} \ \partial (Id_{E^*} \otimes \mu)$ is an isomorphism from $E^* \otimes E^*$ to itself.
\end{prop}
\begin{proof}
\begin{eqnarray*}
\mathbf{P}(cd \otimes ab) = ad \otimes cb + cb \otimes ad -(m+r) cd \otimes ab.
\end{eqnarray*}
Therefore
\begin{eqnarray*}
\frac{-1}{m+r} \mathbf{P}(cd \wedge ab) = cd \wedge ab,
\end{eqnarray*}
and since $m+r > 2$,
\begin{eqnarray*}
\frac{2}{4 - (m+r)^2} \big( \mathbf{P}(ad \odot cb) + \frac{m+r}{2} \mathbf{P}(cd \odot ab) \big) = cd \odot ab.
\end{eqnarray*}
showing that $\mathbf{P}$ is surjective, and, equivalently, bijective.
\end{proof}

All this implies that the $\partial ( V^* \otimes \mathfrak{g}^{(1)} )$ component of the curvature is of Ricci-type. Then the whole curvature must be of Ricci-type, except for the minimal Segre algebras. We will deal with those in the next chapter.

\section{Minimal Segre algebras}
There is no uniform terminology for algebras of this type. As the general algebras $\mathbb{C} \oplus \mathfrak{sl}(m, \mathbb{C}) \oplus \mathfrak{sl}(r,\mathbb{C})$ are sometimes called Segre structures, I have elected to call them `minimal Segre' when $r$ is minimal -- though they are sometimes referred to as `paraconformal'. Recall that these are the algebras
\begin{eqnarray*}
\mathbb{C} \oplus \mathfrak{sl}(m_1, \mathbb{C}) \oplus \mathfrak{sl}(2,\mathbb{C}), \\
\mathbb{R} \oplus \mathfrak{sl}(m_2, \mathbb{R}) \oplus \mathfrak{sl}(2,\mathbb{R}), \\
\mathbb{R} \oplus \mathfrak{sl}(m_3, \mathbb{H}) \oplus \mathfrak{sl}(1,\mathbb{H}).
\end{eqnarray*}
Notice that for $m_1 = m_2$ and $m_1 = 2 m_3$, the second two algebras are real forms of the first. Furthermore, Ricci-flatness forces the preservation of a complex volume-form by Lemma \ref{complex:ricci}; we shall consequently only have to use the complex algebra $\mathfrak{g} = \mathfrak{sl}(m, \mathbb{C}) \oplus \mathfrak{sl}(2,\mathbb{C})$ in this section.

\begin{theo} \label{last:theorem}
Let $\nabla$ be a Ricci-flat affine connection whose holonomy is contained in $\mathfrak{sl}(m, \mathbb{C}) \oplus \mathfrak{sl}(2,\mathbb{C})$. Then its holonomy is contained in $\mathfrak{sl}(m, \mathbb{C})$.
\end{theo}
As a direct consequence of this theorem, we can say that any subalgebra or real form of $\mathfrak{g}$ acting irreducibly, cannot be a Ricci-flat holonomy algebra. This concerns the following algebras:
\begin{eqnarray*}
\begin{array}{|c|c|}
\hline
\hline
\mathrm{Algebra} & \mathrm{Representation} \\
\hline
\hline
& \\
\mathfrak{sl}(2, \mathbb{C}) \oplus \mathfrak{sl}(n,\mathbb{C}) & \mathbb{C}^2 \otimes_{\mathbb{C}} \mathbb{C}^n \cong \mathbb{C}^{2n}, \ \ n \geq 3\\
& \\
\mathfrak{sl}(2, \mathbb{R}) \oplus \mathfrak{sl}(n,\mathbb{R}) & \mathbb{R}^2 \otimes_{\mathbb{R}} \mathbb{R}^n \cong \mathbb{R}^{2n}, \ \ n \geq 3 \\
& \\
\mathfrak{sl}(1, \mathbb{H}) \oplus \mathfrak{sl}(n,\mathbb{H}) & \mathbb{H}^1 \otimes_{\mathbb{H}} \mathbb{H}^n \cong \mathbb{R}^{4n}, \ \ n \geq 2 \\
& \\
\mathfrak{sl}(2, \mathbb{C}) \oplus \mathfrak{sp}(2n,\mathbb{C}) & \mathbb{C}^2 \otimes_{\mathbb{C}} \mathbb{C}^{2n} \cong \mathbb{C}^{4n}, \ \ n \geq 2 \\
& \\
\mathfrak{sl}(2, \mathbb{R}) \oplus \mathfrak{sp}(2n,\mathbb{R}) & \mathbb{R}^2 \otimes_{\mathbb{R}} \mathbb{R}^{2n} \cong \mathbb{R}^{4n}, \ \ n \geq 2 \\
&\\
\mathfrak{sp}(1) \oplus \mathfrak{sp}(p,q) & \mathbb{H} \otimes_{\mathbb{H}} \mathbb{H}^{(p,q)} \cong \mathbb{R}^{(4q,4q)}, \ \ p+q \geq 2 \\
&\\
\hline
\end{array}
\end{eqnarray*}
Though of course in that last case the result -- that a Ricci-flat quaternionic-K\"ahler manifold is hyper-K\"ahler -- is well known, \cite{4KM} and \cite{old4KM}.

In order to prove this theorem, we shall use the quaternionic approach from paper \cite{quaternionic}, modified to incorperate the full complex case.

Let $\mathcal{G}$ be the frame bundle for the $\mathfrak{g}$-structure, and let $J_1$, $J_2$ and $J_3$ be sections of the bundle $Q = \mathcal{G} \times_G \mathfrak{sl}(2, \mathbb{C})$, chosen so they obey the quaternionic identities $J_{\alpha} J_{\beta} = - \delta_{\alpha \beta} Id + \epsilon_{\alpha \beta \gamma} J_{\gamma}$. Thus the complex span of these elements cover all of $Q$.

Let $\nabla$ be any connection associated to this $\mathfrak{g}$-structure. The curvature $R^{\nabla}$ of $\nabla$ decomposes as
\begin{eqnarray*}
R^{\nabla '} + \Omega^1 J_1 + \Omega^2 J_2 + \Omega^3 J_3,
\end{eqnarray*}
where $R^{\nabla '}$ is a curvature terms with values in $\mathcal{G} \times_G \mathfrak{sl}(m, \mathbb{C})$, and the $\Omega_{\alpha}$ are sections of $\wedge^2 T^* \otimes \mathbb{C}$.

Note the commutator relation
\begin{eqnarray*}
\big[ R^{\nabla}, J_{\alpha} \big] = 2 \left( \Omega^{\gamma} J_{\beta} - \Omega^{\beta} J_{\gamma} \right)
\end{eqnarray*}
where $(\alpha, \beta, \gamma)$ is a cyclic permutation of $(1,2,3)$. Let $\Omega^{\alpha'}$ and $\Omega^{\alpha''}$ be the real and imaginary parts of $\Omega^{\alpha}$. Since all elements of $\mathfrak{g}$ are trace-free, we may calculate the $\Omega^{\alpha}$ using the formula
\begin{eqnarray*}
\Omega^{\alpha'}_{(X,Y)} = - \frac{1}{4m} Tr(R^{\nabla}_{(X,Y)} \circ J_{\alpha})
\end{eqnarray*}
and
\begin{eqnarray*}
\Omega^{\alpha''}_{(X,Y)} = \frac{1}{4m} Tr(R^{\nabla}_{(X,Y)} \circ i J_{\alpha}).
\end{eqnarray*}
Note that these traces are real traces.

There are two other operators we shall be needing: the $i$-linearity operator  $\widetilde{}$  and the operator  $\widehat{}$ , the hermitian operator with respect to the complex structure. In details, for any section $F$ of $\wedge^2 T^* \otimes \mathbb{C}$,
\begin{eqnarray*}
\widetilde{F} (X,Y) = \frac{1}{2} (F (X,Y) - i F (X,iY))
\end{eqnarray*}
while
\begin{eqnarray*}
\widehat{F} (X,Y) = \frac{1}{4} \big( F (X,Y) - \sum_{k=1}^3 F (J_{k} X, J_k Y) \big).
\end{eqnarray*}
It is easy to see that both these operators are projections, i.e.~square to themselves. The $i$-linearity operator has certain interesting properties; indeed
\begin{lemm} \label{last:lemma}
If $F$ is a section of $\wedge^{(2,0)} T^* \otimes \mathbb{C}$ -- the tensor product in this expression is real -- then $\widetilde{F}$ is skew-symmetric. If $F$ is a section of $\wedge^{(1,1)} T^* \otimes \mathbb{C}$, then
\begin{eqnarray*}
F(X,Y) = \widetilde{F}(X,Y) - \widetilde{F}(Y,X).
\end{eqnarray*}
\end{lemm}
\begin{proof}
\begin{eqnarray*}
\widetilde{F}(X,Y) + \widetilde{F}(Y,X) &=& \frac{1}{2} \big( F(X,Y) + F(Y,X) - i F(X,iY) - iF(Y,iX) \big) \\ &=& 0 + \frac{-i}{2}\big( F(X,iY) + F(iY,X) \big) \\
&=& 0,
\end{eqnarray*}
if $F \in \Gamma(\wedge^{(2,0)} T^* \otimes \mathbb{C})$. On the other hand,
\begin{eqnarray*}
\widetilde{F}(X,Y) - \widetilde{F}(Y,X) &=& \frac{1}{2} \big( F(X,Y) - F(Y,X) - i F(X,iY) + iF(Y,iX) \big) \\ &=& F(X,Y) + \frac{-i}{2}\big( F(X,iY) + F(iY,X) \big) \\
&=& F(X,Y),
\end{eqnarray*}
if $F \in \Gamma(\wedge^{(1,1)} T^* \otimes \mathbb{C})$.
\end{proof}
\noindent We now aim to show that
\begin{prop}
If $\nabla$ is Ricci-flat, then $\Omega^{\alpha} = 0$ for all $\alpha$.
\end{prop}
\begin{proof}
We shall prove this statement purely algebraically. Since we may, as in section \ref{complex:algebras}, split the curvature module into two components, the $i$-symmetric and $i$-hermitian components, both obeying the Bianchi identity and whose Ricci tensors are respectively $i$-symmetric and $i$-hermitian, it suffices to prove this result in the two cases where $\nabla$ is assumed to have purely $i$-symmetric and purely $i$-hermitian curvature.

We shall deal with the first case first. Notice that this implies that $\Omega^{\alpha}$ is a section of $\wedge^{(2,0)} T^* \otimes \mathbb{C}$.

Let $(E^k)$ be a local frame on the manifold, with dual frame $(e_k)$. Then using the Bianchi identity, the function $-4m \Omega^{\alpha'}_{(J_{\alpha} X, J_{\alpha} Y)}$ is equal to
\begin{eqnarray*}
Tr(R^{\nabla}_{(J_{\alpha} X, J_{\alpha} Y)} \circ J_{\alpha}) &=& \sum_k \left( R^{\nabla}_{(J_{\alpha} X, J_{\alpha} Y)} J_{\alpha} E^k \right) \llcorner e_k \\
&=& - \sum_k \left( R^{\nabla}_{(J_{\alpha} Y, J_{\alpha} E^k)} J_{\alpha} X \right) \llcorner e_k - \sum_k \left( R^{\nabla}_{(J_{\alpha} E^k, J_{\alpha} X)} J_{\alpha} Y \right) \llcorner e_k \\
&=& - \sum_k \left( R^{\nabla}_{(J_{\alpha} Y, J_{\alpha} E^k)} X \right) \llcorner J_{\alpha} e_k - \sum_k \left( R^{\nabla}_{(J_{\alpha} E^k, J_{\alpha} X)} Y \right) \llcorner J_{\alpha} e_k \\
&& - \sum_k \left( 2\Omega^{\gamma}_{(J_{\alpha}Y, J_{\alpha} E^k)} J_{\beta} X - 2\Omega^{\beta}_{(J_{\alpha}Y, J_{\alpha} E^k)} J_{\gamma} X \right) \llcorner e_k \\
&& - \sum_k \left( 2\Omega^{\gamma}_{(J_{\alpha} E^k, J_{\alpha} X)} J_{\beta} Y - 2\Omega^{\beta}_{(J_{\alpha}E^k, J_{\alpha} X)} J_{\gamma} Y \right) \llcorner e_k \\
&=& - \mathsf{Ric}(J_{\alpha}Y, X) + \mathsf{Ric}(J_{\alpha}X, Y) \\
&& - 2 \big( \Omega^{\gamma'}_{(J_{\alpha}Y, J_{\gamma}X)} + \Omega^{\beta'}_{(J_{\alpha}Y, J_{\beta}X)} + \Omega^{\gamma'}_{(J_{\gamma}Y, J_{\alpha}X)} + \Omega^{\beta'}_{(J_{\beta}Y, J_{\alpha}X)} \big) \\
&& - 2 \big( \Omega^{\gamma''}_{(J_{\alpha}Y, i J_{\gamma}X)} + \Omega^{\beta''}_{(J_{\alpha}Y, i J_{\beta}X)} + \Omega^{\gamma''}_{(iJ_{\gamma}Y, J_{\alpha}X)} + \Omega^{\beta''}_{(i J_{\beta}Y, J_{\alpha}X)} \big)
\end{eqnarray*}
The $\mathsf{Ric}$ terms disappear, of course, and using the corresponding expression for $-4m \Omega^{\alpha''}_{(J_{\alpha} X, J_{\alpha} Y)}$, one gets the equation
\begin{eqnarray} \label{omega:symmetric}
-4m \Omega^{\alpha}_{(J_{\alpha} X, J_{\alpha} Y)} &=& - 4 \big( \widetilde{\Omega}^{\gamma}_{(J_{\alpha}Y, J_{\gamma}X)} + \widetilde{\Omega}^{\beta}_{(J_{\alpha}Y, J_{\beta}X)} + \widetilde{\Omega}^{\gamma}_{(J_{\gamma}Y, J_{\alpha}X)} + \widetilde{\Omega}^{\beta}_{(J_{\beta}Y, J_{\alpha}X)} \big),
\end{eqnarray}
since the $\Omega^{\alpha}$ are $i$-symmetric. Notice that this equation implies that $\Omega^{\alpha}$ is completely $i$-symmetric, i.e.~that $\Omega^{\alpha} = \widetilde{\Omega}^{\alpha}$. By replacing $Y$ with $J_{\alpha} Y$ and defining $\Omega = \sum_k \Omega^k_{(\cdot, J_k \cdot)}$, we may rewrite this equation as
\begin{eqnarray} \label{alpha:equation}
(m-2) \Omega^{\alpha}_{Y, (J_{\alpha} X)} + \Omega_{(J_{\alpha} X, J_{\alpha}Y)} + \Omega_{(Y,X)} = 0.
\end{eqnarray}
By summing over $\alpha = 1, 2, 3$, we get the identity
\begin{eqnarray*}
(m+1) \Omega_{(Y, X)} + \sum_{\alpha} \Omega_{(J_{\alpha}X, J_{\alpha} Y)} =0,
\end{eqnarray*}
from which it follows that, if $\Omega^s$ and $\Omega^a$ are the symmetric and anti-symmetric parts of $\Omega$,
\begin{eqnarray*}
-m \Omega^s = 4 \widehat{\Omega}^s \ \ \mathrm{and} \ \ (m+2) \Omega^a = 4 \widehat{\Omega}^a.
\end{eqnarray*}
which, since $m + 2 \neq 4$ and $-m \neq 4$, implies that $\Omega = 0$ and hence, by Equation (\ref{alpha:equation}), that $\Omega^{\alpha} = 0$.

We now turn to the $i$-hermitian piece, for which the proof starts in the same manner, except that Equation (\ref{omega:symmetric}) becomes
\begin{eqnarray*}
-4m \Omega^{\alpha}_{(J_{\alpha} X, J_{\alpha} Y)} &=& - 4 \big( \widetilde{\Omega}^{\gamma}_{(J_{\alpha}Y, J_{\gamma}X)} + \widetilde{\Omega}^{\beta}_{(J_{\alpha}Y, J_{\beta}X)} - \widetilde{\Omega}^{\gamma}_{(J_{\alpha}X, J_{\gamma}Y)} - \widetilde{\Omega}^{\beta}_{(J_{\alpha}X, J_{\beta}Y)} \big).
\end{eqnarray*}
Notice the exchange of indices and signs in the last two terms. Since the section $\widetilde{F}$ is defined by the relation $\widetilde{F}(X,iY) = i \widetilde{F}(X,Y)$ and $\widetilde{F}$ remains $i$-hermitian if $F$ is, we may deduce that
\begin{eqnarray*}
m \widetilde{\Omega}^{\alpha}_{(J_{\alpha} X, J_{\alpha} Y)} &=& - \widetilde{\Omega}^{\gamma}_{(J_{\alpha}X, J_{\gamma}Y)} - \widetilde{\Omega}^{\beta}_{(J_{\alpha}X, J_{\beta}Y)},
\end{eqnarray*}
or, equivalently, after replacing $X$ with $J_{\alpha} X$,
\begin{eqnarray*}
(m-1) \widetilde{\Omega}^{\alpha}_{(X, J_{\alpha} Y)} + \widetilde{\Omega}_{(X, Y)} = 0.
\end{eqnarray*}
Summing over $\alpha$ gives us $(m+2) \widetilde{\Omega} = 0$ and, consequently,
\begin{eqnarray*}
\widetilde{\Omega}^{\alpha} = 0.
\end{eqnarray*}
And then the relation $\Omega^{\alpha} (X,Y) = \widetilde{\Omega}^{\alpha}(X,Y) - \widetilde{\Omega}^{\alpha}(Y,X)$ from Lemma \ref{last:lemma} gives us the required vanishing of $\Omega^{\alpha}$.
\end{proof}
And this is all we need to prove Theorem \ref{last:theorem}.

\section{The case of $E_6$}

Also present in the table of possible irreducible torsion-free affine holonomy algebras are various sub-algebras and real forms of
\begin{eqnarray*}
\mathfrak{g} = \mathbb{C}^* \cdot \mathfrak{e}_6^{\mathbb{C}}.
\end{eqnarray*}
We aim to prove that $\mathfrak{g}$ is of Ricci-type, and that consequently all sub-algebras and real forms of it are. The representation space of $\mathfrak{g}$ is
\begin{eqnarray*}
V \cong \mathbb{C}^{27}.
\end{eqnarray*}
This is the standard representation space of $\mathfrak{g}$. The algebra $\mathfrak{e}_6^{\mathbb{C}}$ is in fact defined as the maximal algebra preserving a certain non-degenerate cubic $\Psi$ on $V$ \cite{E6}. Non-degeneracy means that the $\Psi$-induced maps, $V \to \odot^2 V^*$ and $\odot^2 V \to V^*$ are of maximum rank. The full algebra $\mathfrak{g}$ must preserve $\Psi$ up to scaling.

The Dynkin diagram of $\mathfrak{e}_6^{\mathbb{C}}$ has six nodes, and the maximal weights are given by sextuplets of non-negative integers. In this optic,
\begin{eqnarray*}
V &=& (1,0,0,0,0,0) \\
V^* &=& (0,0,0,0,1,0).
\end{eqnarray*}
The dual representation of $\mathfrak{e}_6^\mathbb{C}$ on $V^*$ must preserve a non-degenerate cubic $\Theta \in \odot^3 V$. We choose the scale of $\Theta$ by requiring, in abstract index notation,
\begin{eqnarray*}
\Psi_{jkl} \Theta^{jkl} = 27.
\end{eqnarray*}
We will use \cite{BoringTable} in order to calculate various tensor products of representations of $\mathfrak{e}_6^\mathbb{C}$. Using $\Theta$, there is a decomposition of
\begin{eqnarray*}
\odot^2 V^* = V \oplus U.
\end{eqnarray*}
Using \cite{BoringTable}, one has that $U$ is irreducible and
\begin{eqnarray*}
U = (0,0,0,0,2,0).
\end{eqnarray*}
This decomposition implies that,
\begin{eqnarray*}
\Psi_{jkm} \Theta^{jkl} = Id_m^l,
\end{eqnarray*}
as a map $V \to V$ or $V^* \to V^*$. Similarly,
\begin{eqnarray*}
\Psi_{jpm} \Theta^{jkl} = \Pi_{pm}^{kl},
\end{eqnarray*}
where $\Pi$ is the projection of $\odot^2 V^*$ onto its submodule $V$, along $U$. Using \cite{BoringTable}, we can decompose $\odot^2 V^* \otimes V$,
\begin{eqnarray*}
\odot^2 V^* \otimes V =&& (1,0,0,0,2,0) \oplus (0,0,0,0,1,1) \oplus (0,0,0,0,2,0) \\
& \oplus & (0,1,0,0,0,0) \oplus U^* \oplus 2V^*.
\end{eqnarray*}
Now, we know by \cite{CHQM} and \cite{CIH} that
\begin{eqnarray*}
\mathfrak{g}^{(1)} &=& V^* \\
H^{1,2}(\mathfrak{g}) &=& 0.
\end{eqnarray*}
Therefore the module $\mathfrak{g}^{(1)}$ is contained inside the two $V^*$ components of the previous decomposition. Then define two maps $V^* \to \odot^2 V^* \otimes V$ by
\begin{eqnarray*}
\mu_1(v_j) &=& \Psi_{kmr} \Theta^{jkl} v_j \\
\mu_2(v_j) &=& v_{j} Id_{k}^l + v_k Id_j^l
\end{eqnarray*}
where $(jk)$ denotes symmetrisation of the indices.
\begin{lemm} \label{mu:prop}
$\mu_1$ and $\mu_2$ are injective and non-isomorphic.
\end{lemm}
\begin{proof}
The traces of $\mu_1$ and $\mu_2$ are
\begin{eqnarray*}
\mathrm{trace} \ \mu_1(v_j) &=& v_j \\
\mathrm{trace} \ \mu_2(v_j) &=& (28) v_j,
\end{eqnarray*}
proving that both are non-zero, hence (as $V^*$ is irreducible) injective. Now contract them with an element $w_l$ of $V^*$:
\begin{eqnarray*}
\mu_1(v_j) \ \llcorner  \ w_l &=& \Psi_{kmr} \Theta^{jkl} v_j w_l = \Pi^{jk}_{mr} (v_j w_l)\\
\mu_2(v_j) \ \llcorner \ w_l &=& v_{j} w_{l} + w_{j} v_{l},
\end{eqnarray*}
and these two elements cannot be isomorphic for general $v_j$ and $w_l$.
\end{proof}
Consequently any map $\nu: V^* \to \odot^2 V^* \otimes V$ is given as
\begin{eqnarray*}
\nu = \lambda_1 \mu_1 + \lambda_2 \mu_2,
\end{eqnarray*}
for complex constants $\lambda_1$ and $\lambda_2$.

\begin{prop}
Let $\nu: V^* \to \mathfrak{g}^{(1)}$ be an invariant isomorphic map. Then $\nu$ has $\lambda_1 = -\lambda_2$.
\end{prop}
\begin{proof}
We shall leave abstract index notation to the side for the moment, and we will need to provide a more explicit description of $\Psi$ and $\Theta$. There is an inclusion $\mathfrak{sp}(8, \mathbb{C}) \subset \mathfrak{e}_6^\mathbb{C}$ described as follows. Let $\omega$ be the symplectic form for $\mathfrak{sp}(8, \mathbb{C})$. Then there is a map
\begin{eqnarray*}
\wedge^2 \mathbb{C}^{8*} \to \wedge^8 \mathbb{C}^{8*},
\end{eqnarray*}
by wedging with $\omega^3$. The kernel of this map is 27-dimensional; we shall call it $V^*$, as it is the dual natural representation space of $\mathfrak{e}_6^\mathbb{C}$. To confirm this, consider the non-degenerate cubic $\Theta$ defined on it by
\begin{eqnarray*}
\Theta(a,b,c) \omega^4 = a \wedge b \wedge c \wedge \omega.
\end{eqnarray*}
This cubic is obviously preserved by $\mathfrak{sp}(8, \mathbb{C})$, giving us the required inclusion. Let $(X^j)$ be a basis for $V$, with dual basis $(\eta_j)$. We may express $\omega$ as
\begin{eqnarray*}
\eta_1 \wedge \eta_2 + \eta_3 \wedge \eta_4 + \eta_5 \wedge \eta_6 + \eta_7 \wedge \eta_8.
\end{eqnarray*}
Consequently a basis for $V$ is given by
\begin{eqnarray*}
& \eta_1 \wedge \eta_2 - \eta_3 \wedge \eta_4, \\
& \eta_1 \wedge \eta_2 - \eta_5 \wedge \eta_6, \\
& \eta_1 \wedge \eta_2 - \eta_7 \wedge \eta_8, \\
& \eta_{\alpha} \wedge \eta_{\beta},
\end{eqnarray*}
where $\alpha$ and $\beta$ are numbers chosen from distinct sets in the collection $\{ 1, 2 \}$, $\{3,4 \}$, $\{5, 6 \}$, $\{7, 8 \}$.

We will work on an explicit example to find the (unique) relation between $\lambda_1$ and $\lambda_2$. So let
\begin{eqnarray*}
a &=& \eta_4 \wedge \eta_6 \\
b &=& \eta_1 \wedge \eta_4 \\
c &=& \eta_3 \wedge \eta_6 \\
d &=& \eta_2 \wedge \eta_5.
\end{eqnarray*}
As a consequence, $\Theta(ab) = \Theta(bc) = 0$ and $\Theta(bc) \neq 0$.

Now
\begin{eqnarray*}
a \wedge d \wedge \omega &=& - \eta_2 \wedge \eta_4 \wedge \eta_5 \wedge \eta_6 \wedge \eta_7 \wedge \eta_8.
\end{eqnarray*}
The only basis element this wedges with in a non-trivial way is $\eta_1 \wedge \eta_3$, to give $1$. Consequently,
\begin{eqnarray*}
\Theta(ad) = X^1 \wedge X^3 = Z.
\end{eqnarray*}
We now aim to calculate $\Pi(ad) = \Psi (Z)$. If $(Y^{\sigma})$ is a basis for $\odot^2 V^*$ with dual basis $(y_{\sigma})$,
\begin{eqnarray*}
\Psi(Z) = \sum_{\sigma} \Psi(Z Y^{\sigma}) y_{\sigma}.
\end{eqnarray*}
If $(Y^{\sigma})$ is the tensor product of the basis elements of $V^*$, the only $Y^{\sigma}$ such that $\Psi(Z Y^{\sigma}) \neq 0$ are
\begin{eqnarray*}
X^2 \wedge X^5 &\odot& X^4 \wedge X^6 \\
X^2 \wedge X^6 &\odot& X^4 \wedge X^5 \\
X^2 \wedge X^7 &\odot& X^4 \wedge X^8 \\
X^2 \wedge X^8 &\odot& X^4 \wedge X^7 \\
X^2 \wedge X^4 &\odot& \frac{1}{2} (X^1 \wedge X^2 - X^5 \wedge X^6) \\
X^2 \wedge X^4 &\odot& \frac{1}{2} (X^1 \wedge X^2 - X^7 \wedge X^8).
\end{eqnarray*}
Consequently there is an $ \eta_2 \wedge \eta_5 \odot \eta_4 \wedge \eta_6$ summand in $\Pi(ad)$. In other words, if $W = X^4 \wedge X^6$,
\begin{eqnarray*}
\Pi(ad) \ \llcorner \ W = \eta_2 \wedge \eta_5 = d.
\end{eqnarray*}
Note that $\Psi ( d, b, c ) = -1$.

An element $e$ of $\mathfrak{g}$ preserves $\Theta$ up to scale. In our case $\Theta(abc) = 0$, so there is no issue of scale. Explicitly,
\begin{eqnarray*}
0 = \Theta ((e \cdot a),b,c) + \Theta (a,(e \cdot b),c) + \Theta (a,b, (e \cdot c))
\end{eqnarray*}
since $\Theta$ is zero on $abc$. Because of the choices of $a$, $b$, and $c$ that we made, this ensured that
\begin{eqnarray*}
0 = \Theta ((e \cdot a),b,c).
\end{eqnarray*}
The above formula must hold replacing $e$ with the element of $\mathfrak{e}_6^{\mathbb{C}}$ that is $\nu(d) \ \llcorner \ W$. This implies
\begin{eqnarray*}
0 &=& \lambda_1 \Theta(d,b,c) + \lambda_2 \Theta(d,b,c) \\
&=& -(\lambda_1 + \lambda_2),
\end{eqnarray*}
since
\begin{eqnarray*}
\mu_1(d) \ \llcorner \ a \ \llcorner \ W = \Pi(ad) \ \llcorner W &=& d \\
\mu_2(d) \ \llcorner \ a \ \llcorner W = W(d) a + W(a) d &=& d.
\end{eqnarray*}
\end{proof}
\begin{theo}
The algebra
\begin{eqnarray*}
\mathfrak{g} = \mathbb{C}^* \cdot \mathfrak{e}_6^{\mathbb{C}}.
\end{eqnarray*}
acting on $V \cong \mathbb{C}^{27}$, is of Ricci-type.
\end{theo}
\begin{proof}
The curvature bundle of $\mathfrak{g}$ is
\begin{eqnarray*}
K(\mathfrak{g}) = \partial ( V^* \otimes m(V^*))
\end{eqnarray*}
As before, define $\mathbf{R}$ as the Ricci-trace map $V^* \otimes V^*$ to itself. Then by the properties of $\mu_1$ and $\mu_2$ expounded in Lemma \ref{mu:prop},
\begin{eqnarray*}
\mathbf{R} (w \otimes v) = \lambda_1 \Pi(w \odot v) -2 \lambda_1 w \odot v + 27 \lambda_1  w \otimes v.
\end{eqnarray*}
Now the image of $\mathbf{R}$ is not symmetric, and $\wedge^2 V^*$ is an irreducible representation of $\mathfrak{e}_6^{\mathbb{C}}$; consequently the entire $\wedge^2 V^*$ is in the image of $\mathbf{R}$. Now looking at the symmetric part:
\begin{eqnarray*}
\mathbf{R} (w \odot v) = \lambda_1 \big( \Pi(w \odot v) + 25 w \odot v \big).
\end{eqnarray*}
Consequently
\begin{eqnarray*}
\mathbf{R} (\Pi(w \odot v)) = 26 \lambda_1 \Pi(w \odot v),
\end{eqnarray*}
and
\begin{eqnarray*}
\mathbf{R} ((1- \Pi) (w \odot v)) = 25 \lambda_1 (1- \Pi) (w \odot v).
\end{eqnarray*}
Consequently, as $\lambda_1 \neq 0$ since $\nu$ is non-trivial, $\mathbf{R}$ is an isomorphism, and
\begin{eqnarray*}
\mathfrak{g} = \mathbb{C}^* \cdot \mathfrak{e}_6^{\mathbb{C}}
\end{eqnarray*}
is of Ricci-type.
\end{proof}

\section{Low dimensional cases} \label{RicciF:holo}
In this section, we aim to finish the classification of which holonomy algebras acting irreducibly can correspond to a Ricci-flat connection. For though we have excluded many holonomy algebras from being Ricci-flat, paper \cite{mepro2} has constructed Ricci-flat cones for most of the others, we have not settle the existence of general Ricci-flat connections in some cases. These are the algebras concerned (those that can be Ricci-flat have been marked with a $*$):
\begin{eqnarray*}
\begin{array}{|c|c|c|}
\hline
\hline
\textrm{algebra }\mathfrak{g} & \textrm{representation V} & \textrm{Dimensions} \\
\hline
\hline
& & \\
\mathfrak{so}(p,q) & \mathbb{R}^{(p,q)} &  p+q = 3,4^*  \\
\mathfrak{so}(n, \mathbb{C}) & \mathbb{C}^n & n = 3,4^* \\
\mathfrak{sl}(n, \mathbb{R}) & \mathbb{R}^n & n= 2 \\
\mathfrak{sl}(n, \mathbb{C}) & \mathbb{C}^n & n=  1,2^* \\
\mathfrak{sl}(n, \mathbb{H}) & \mathbb{H}^n & n=1^* \\
& & \\
\hline
\end{array}
\end{eqnarray*}
Of these, we can immediately exclude $\mathfrak{sl}(2, \mathbb{R})$ and $\mathfrak{sl}(1, \mathbb{C})$, as any Ricci-flat two-manifold is flat. In contrast, any $\mathfrak{sl}(n, \mathbb{H})$ connection must be Ricci-flat by definition. Manifolds with holonomy $\mathfrak{so}(p,q)$, $p+q = 3$, have vanishing Weyl tensor as all three-manifolds do. However, a Ricci-flat manifold has full curvature contained in the Weyl tensor. Thus Ricci-flat manifolds with these holonomies must be flat. The result holds, similarly, in the holomorphic category of $\mathfrak{so}(3, \mathbb{C})$.

For the case of $\mathfrak{sl}(2, \mathbb{C})$, let $x$ and $y$ be complex coordinates with corresponding holomorphic vector fields $X, Y \in \Gamma(T_{\mathbb{C}})$. Then define the connection $\nabla$ as
\begin{eqnarray*}
\nabla_{\overline{X}} \overline{X} = \overline{\nabla_{\overline{X}}X } = \overline{\nabla_{X} \overline{X}} = \nabla_X X &=& fY
\end{eqnarray*}
and all other terms involving $X, Y$ and their conjugates are zero. Here $f$ is a complex-valued function that is independent of $y$ (i.e.~$Yf = \overline{Y}f = 0$). This $\nabla$ is a torsion-free connection respecting the real structure on $T_{\mathbb{C}}$ -- consequently equivalent to a real connection representing the corresponding complex structure on $T$. The curvature of $\nabla$ is given by
\begin{eqnarray*}
R_{\overline{X} X} X &=& (\overline{X}f)Y - (Xf) \overline{Y},
\end{eqnarray*}
the corresponding $R_{\overline{X} X} \overline{X}$ term, and all other curvature terms are zero. This makes $\nabla$ Ricci-flat. We now use $f$ as a bump function to smoothly move $\nabla$ to a flat connection (moving along the $x$ direction, of course), while remaining Ricci-flat and complex along the way. Given two copies of this manifold, we may glue two distinct patches on each manifold while identifying the coordinate $x$ with $y$ and $y$ with $ix$. This may be done so that the resulting structure is a manifold $M$. We thus have a complex, Ricci-flat connection $\nabla$ such that, directly from the curvature, we have the holonomy elements
\begin{eqnarray*}
X \to Y,
\end{eqnarray*}
\noindent and
\begin{eqnarray*}
Y \to iX.
\end{eqnarray*}
And these two elements generate the full $\mathfrak{sl}(T, \mathbb{C})$ holonomy.

In order to generate the remaining holonomy algebras, we turn to the Schwarzschild metric \cite{Schmetric}. In this (Lorentzian) case, the metric is
\be
g = - C dt^2 +\frac{1}{C} dr^2 + r^2(d \theta^2 + \sin{\theta} d \psi^2).
\ee
Where $C = 1 - \frac{2M}{r}$ for some mass $M$. There is also an Euclidean Schwarzschild metric (given by replacing $-dt^2$ with $dt^2$) and a split Schwarzschild metric (given by replacing $d \psi^2$ with $-d \psi^2$). Of course, since the metric is real-analytic in the coordinates, there is also a complex Schwarzschild metric, my considering $t, r, \theta $ and $\psi$ as complex coordinates.

Then elementary but laborious calculations establish that all these metrics are Ricci-flat, and all have maximal holonomy -- it turns out that the curvature tensor is enough to generate the full holonomy algebra.

\end{document}